\scrollmode
\documentclass[12pt]{amsart}
\usepackage{verbatim}
\usepackage{eucal}
\usepackage{amssymb}
\usepackage{mathrsfs}
\usepackage{graphicx}
\usepackage{psfrag}
\usepackage{cite}

\newif \ifwide
\newif \ifavnermargin
\def \makemargins{
\ifwide
        \oddsidemargin .25in
        \evensidemargin .25in
        \textwidth 6.00in
\else
\fi
\ifavnermargin
        \headheight=7pt
        \textheight=574pt
        \textwidth=432pt
        \topmargin=14pt
        \oddsidemargin=18pt
        \evensidemargin=18pt
\else
\fi
}

\avnermargintrue
\makemargins

\theoremstyle{plain}
\newtheorem{theorem}[subsection]{Theorem}
\newtheorem{proposition}[subsection]{Proposition}
\newtheorem{lemma}[subsection]{Lemma}

\theoremstyle{definition}
\newtheorem{example}[subsection]{Example}

\theoremstyle{remark}
\newtheorem{remark}[subsection]{Remark}

\newcount\timehh\newcount\timemm
\timehh=\time
\divide\timehh by 60 \timemm=\time
\newcommand{\draftauthor}[1]{\author{#1
    {
      --- \protect \protect\sc\today\ ---
      \ifnum\timehh<10 0\fi\number\timehh\,:\,\ifnum\timemm<10 0\fi\number\timemm
      \protect \, \, \protect 
    }
  }
}
\newcommand{\RR}{{\mathbb R}}

\newcommand{\ZZ}{{\mathbb Z}}
\newcommand{\QQ}{{\mathbb Q}}

\newcommand{\FF}{{\mathbb F}}
\newcommand{\Proj}{{\mathbb P}}

\newcommand{\fH}{{\mathfrak{H}}}

\newcommand{\PPP}{{\mathscr{P}}}

\newcommand{\FFF}{{\mathscr{F}}}
\newcommand{\TTT}{{\mathscr{T}}}
\newcommand{\HHH}{{\mathscr{H}}}

\newcommand{\MMM}{{\mathscr{M}}}

\newcommand{\LLL}{{\mathscr{L}}}

\renewcommand{\tilde}{\widetilde}

\newcommand{\one}{{\mathbf{1}}}
\newcommand{\dee}[1]{\partial /\partial #1}
\newcommand{\fractwo}[2]{(#1)/ (#2)}
\newcommand{\pint}[1]{[#1]_{p}}
\newcommand{\nonsing}[1]{{#1}_{\text{ns}}}
\newcommand{\sing}[1]{{#1}_{\text{s}}}
\newcommand{\ol}[1]{\overline{#1}}
\newcommand{\avg}[1]{{#1}_{\text{avg}}}

\DeclareMathOperator{\Gr}{Gr}
\DeclareMathOperator{\GL}{GL}
\DeclareMathOperator{\SL}{SL}

\DeclareMathOperator{\Vol}{Vol}

\DeclareMathOperator{\Dim}{dim}

\DeclareMathOperator{\Ind}{Ind}

\DeclareMathOperator{\Td}{Td}

\begin{document}

\newif \ifdraft
\def \makeauthor{
\ifdraft
        \draftauthor{PEG and FRV}
\title[Hecke operators and the Ehrhart polynomial (\RCSRevision)]{Lattice polytopes,
Hecke operators, and the Ehrhart polynomial (\RCSRevision)}

\else
\title[Hecke operators and the Ehrhart polynomial]{Lattice polytopes,
Hecke operators, and the Ehrhart polynomial}

\author{Paul E. Gunnells}
\address{Department of Mathematics and Statistics\\
University of Massachusetts\\
Amherst, MA  01003}
\email{gunnells@math.umass.edu}

\author{Fernando Rodriguez Villegas}
\address{Department of Mathematics\\
University of Texas\\
Austin, TX 78712}
\email{villegas@math.utexas.edu}

\fi
}

\drafttrue
\draftfalse
\makeauthor

\ifdraft
        \date{\today}
\else
        \date{May 29, 2004}
\fi


\begin{abstract}
Let $P$ be a simple lattice polytope.  We define an action of the
Hecke operators on $E (P)$, the Ehrhart polynomial of $P$, and describe
their effect on the coefficients of $E (P)$.  We also describe how the
Brion-Vergne formula transforms under the Hecke operators for
nonsingular lattice polytopes $P$.
\end{abstract}
\maketitle

\section{Introduction}\label{introduction}

\subsection{}
Let $L$ be a rank $n$ lattice, embedded in a real $n$-dimensional
vector space $V$.  Let $\PPP (L)$ be the set of $n$-dimensional convex
polytopes in $V$ with vertices in $L$.  For any $P \in \PPP (L)$, and
for any nonnegative integer $t$, let $tP$ be $P$ scaled by the factor
$t$.  Then by a result of Ehrhart \cite{ehrhart}, the function $t\mapsto
\#(tP\cap L)$ is a degree $n$ polynomial with rational coefficients,
called the \emph{Ehrhart polynomial} of $P$.  Hence one can think of
the Ehrhart polynomial as giving a map $E$ from $\PPP (L)$ to the
polynomial ring $\QQ [t]$.

Write $E (P) = \sum_{l=0}^{n} c_{l}t^{l}$.  Formulas for the
coefficients $c_{l}$, in various settings and with varying degrees of
generality, have been given by several authors
\cite{kp,bv1,bv2,diaz-robins,kk,jppp,guillemin,cs}.  Some coefficients
are easy to understand, for example 
\begin{equation}\label{eq:easy.cs}
c_{0}=1, \quad c_{n}=\Vol P, \quad \text{and}\quad 
c_{n-1} = \Vol (\partial P)/2.
\end{equation}
Here $\Vol P$ is taken with respect to the measure that gives a
fundamental domain of $L$ volume 1; if a polytope has dimension less
than $n$, we compute its volume with respect to the lattice obtained
by intersecting its affine hull with $L$.  For a general lattice
polytope, expressions for the Ehrhart coefficients involve not
only volumes but also subtle arithmetic information, namely
\emph{higher-dimensional Dedekind sums} as studied by Carlitz and
Zagier \cite{carlitz, zag-dede}.

\subsection{}
The Ehrhart polynomial depends not just on the combinatorial type of
$P$, but rather on the pair $(P,L)$.  Hence it is natural to consider
how $E (P)$ changes as $L$ is varied.  The theory of automorphic forms
provides a powerful machine to accomplish this, namely the technique
of \emph{Hecke operators}.

Thus let $p$ be a prime, and let $k\leq n$ be a positive integer.
Given a lattice polytope $P$ with Ehrhart polynomial $E (P)$, we
define a new polynomial $T(p,k)E (P)$ as follows.  Let $p^{-1}L$ be
the canonical superlattice of $L$ of coindex $p^{n}$.  We have
$p^{-1}L/L \simeq \FF_{p}^{n}$, and any lattice $M$ satisfying $p^{-1}L
\supsetneq M \supsetneq L$ determines a subspace $\overline M \subset
\FF_{p}^{n}$.  Let $\LLL_{k}$ be the set of such lattices with $\Dim
\overline M = k$.  Then we define
\begin{equation}\label{eq:heckeact}
T (p,k)E (P) = \sum_{{M\in \LLL_{k}}} E (P_{M}),
\end{equation}
where $P_{M}\in \PPP (M)$ denotes the lattice polytope with vertices
in $M$ canonically determined by $P$.  

\subsection{}\label{ss:nu.def}
In this paper we consider the relationship between $T (p,k)E (P)$ and
$E (P)$.  To state our results, we require more notation.  For any
nonnegative integer $l\leq n$, choose and fix an $l$-dimensional
subspace $U$ of $\FF_p^{n}$, and define
\begin{equation}\label{eq:nu.formula}
\nu_{n,k,l}(p) = \sum_{\substack{W\subset \FF_{p}^{n}\\
\Dim W = k}} p^{\Dim W\cap U}.  
\end{equation}
Note that this value is independent of the choice of $U$.
Finally for any polynomial $f\in \QQ [t]$ let $c_{l} (f)$ be the
coefficient of $t^{l}$ in $f$.  Then our first result can be stated as
follows:

\begin{theorem}\label{thm:main1}
We have 
\begin{equation}\label{eq:thm}
c_{l} (T (p,k)E (P))/c_{l} (E (P)) = \nu_{n,k,l}(p),
\end{equation}
independently of $P$.  The ratios $\nu$ satisfy 
\[
\nu_{n,k,l}(p)/\nu_{n,n-k,n-l}(p) = p^{k+l-n}.
\]
Moreover, for each triple $(n,k,l)$, there is a polynomial with
positive coefficients 
\begin{equation}\label{eq:poly.function}
\Phi_{n,k,l} (t)\in \ZZ [t],
\end{equation}
independent of $p$, 
such that $\Phi_{n,k,l} (p) = \nu_{n,k,l} (p)$.
\end{theorem}

The sum \eqref{eq:nu.formula} can be viewed as a sum of $p$-powers
over a certain geometrically-defined stratification of the finite
Grassmannian $\Gr (k,n) (\FF_{p})$, and thus it is not surprising that
for any given $p$ the quantity $\nu_{n,k,l} (p)$ can be expressed as a
integral polynomial in $p$.  However, the existence of $\Phi$, as well
as the statement that it has positive coefficients, does not follow
immediately from \eqref{eq:nu.formula} since the number of terms in
the sum grows with $p$ and since the strata are only locally closed.

As an example of Theorem \ref{thm:main1}, if $l=0$, then $c_{0} (E
(P)) = 1$ for any $P$.  Hence the ratio on the left of \eqref{eq:thm}
is the number of terms in \eqref{eq:heckeact}.  It is well known that
this is the cardinality of $\Gr(k,n) (\FF_{p})$ (cf. Lemma
\ref{lem:cosets}), which equals $\nu_{n,k,0} (p)$.  For further
examples, Table \ref{tab:ev} shows the Hecke eigenvalues that arise
for the Ehrhart coefficients of $4$-dimensional polytopes.

\begin{table}[htb]
\begin{center}
\begin{tabular}{c||c|c|c}
&$T (p,1)$&$T (p,2)$&$T (p,3)$\\
\hline\hline
$c_{4}$&${p^{4} + p^{3} + p^{2} + p}$&${p^{6} + p^{5} + {2}p^{4} + p^{3} + p^{2}}$&${p^{6} + p^{5} + p^{4} + p^{3}}$\\
$c_{3}$&${{2}p^{3} + p^{2} + p}$&${p^{5} + {2}p^{4} + {2}p^{3} + p^{2}}$&${p^{5} + p^{4} + {2}p^{3}}$\\
$c_{2}$&${p^{3} + {2}p^{2} + p}$&${{2}p^{4} + {2}p^{3} + {2}p^{2}}$&${p^{4} + {2}p^{3} + p^{2}}$\\
$c_{1}$&${p^{3} + p^{2} + {2}p}$&${p^{4} + {2}p^{3} + {2}p^{2} + p}$&${{2}p^{3} + p^{2} + p}$\\
$c_{0}$&${p^{3} + p^{2} + p + 1}$&${p^{4} + p^{3} + {2}p^{2} + p + 1}$&${p^{3} + p^{2} + p + 1}$
\end{tabular}
\end{center}
\bigskip
\caption{\label{tab:ev}Eigenvalues for $n=4$.}
\end{table}

\subsection{}
A geometric interpretation of the eigenvalue \eqref{eq:nu.formula} is
the following.  Consider the map
\[
\Vol_{l} \colon \PPP (L) \longrightarrow \QQ
\] 
taking $P$ to the sum of the volumes of all faces of dimension $l$.
Then we can define an action of the Hecke operators on $\Vol_{l}$ as
in \eqref{eq:heckeact}, and one can show that $T(p,k)\Vol_{l} = \nu_{n,k,l}(p)
\Vol_{l}$ (Proposition \ref{prop:vol.tfm}).  Hence Theorem
\ref{thm:main1} says that the $l$th coefficient of the Ehrhart
polynomial transforms under the Hecke operators exactly as the volumes
of the $l$-dimensional faces do.  For another interpretation, in terms
of counting the number of $\FF_{p}$-points on certain varieties, see
Remark \ref{rem:tom}.

\subsection{}\label{ss:main.result}
Recall that an $n$-dimensional lattice polytope is called
\emph{simple} if every vertex meets exactly $n$ edges, and is called
\emph{nonsingular} if for any vertex $v$, the primitive lattice
vectors parallel to the edges emanating from $v$ form a $\ZZ$-basis of
$L$.  Our next result concerns how the Hecke operators interact with
certain formulas for the coefficients of the Ehrhart polynomial in the
special case that $P$ is simple.

Let $\FFF (n-1)$ be the set of facets of $P$, and let $h =
(h_{F})_{F\in \FFF (n-1)}$ be a real multivariable indexed by the
facets of $P$.  Let $P (h)$ be the convex region obtained by parallel
translation of the facets of $P$ by the parameter $h$, normalized by
$P (0) = P$ (\S\ref{ss:fourone}).  For small $h$ the region $P (h)$ is
bounded, and the volume $\Vol P (h)$ is a polynomial function of $h$.

Let $\Sigma$ be the normal fan to $P$ (\S\ref{ss:nor.cones}).  Then
the polytope $P$ determines a differential operator $\Td (\Sigma,
\dee{h})$, called the \emph{Todd operator} (\S\ref{ss:fourseven}).
In the special case that $P$ is nonsingular, this operator is defined
as follows.  Let $\Td (x)$ be the power series expansion of $x/
(1-e^{-x})$, i.e.
\[
\Td (x) = \sum_{j=0}^{\infty}\frac{B_{j}}{j!}x^{j},
\]
where $B_{j}$ are the Bernoulli numbers.  For each $h_{F}$ let $\Td
(\dee{h_{F}})$ be the differential operator obtained by formally
replacing $x$ with $\dee{h_{F}}$ in $\Td (x)$.  Then $\Td (\Sigma ,
\dee{h})$ is defined to be the product
\begin{equation}\label{eq:todd.op.intro}
\Td (\Sigma, \dee{h}) = \prod_{F\in \FFF (n-1)}\Td
(\dee{h_{F}}).
\end{equation}
Note that product may be taken in any order, since the derivatives
mutually commute.  This is an infinite-degree differential operator,
and we denote by
\[
\Td_{l} (\Sigma ,\dee{h})
\]
the homogeneous terms of degree $l$.
By Khovanskii-Pukhlikov \cite{kp} one has 
\[
c_{n-l} (E (P)) = \Td_{l} (\Sigma,\dee{h} ) \Vol P (h) \bigr|_{h=0}.
\]
If the polytope $P$ is simple and not nonsingular, then one must
enlarge \eqref{eq:todd.op.intro} with additional terms involving
higher-dimensional Dedekind sums; the corresponding formula is due to
Brion-Vergne \cite{bv1}.

\subsection{}\label{ss:oneseven} Let $f$ be a face of $P$ of
codimension $\leq l$, and let $\pi = (\pi (F))_{F\supset f}$ be a an
ordered partition of $l$ into positive parts indexed by the facets
containing $f$.  The pair $(f,\pi)$ determines a differential operator
\[
\partial_{f}^{\pi} = \prod_{F\supset f}
(\dee{h_{F}})^{\pi (F)},
\]
and we can 
collect
common terms in \eqref{eq:todd.op.intro} to write
\begin{equation}\label{eq:intro.decomp}
\Td_{l} (\Sigma,\dee{h}) = \sum_{(f,\pi)} A{(f,\pi)}
\partial_{f}^{\pi}.
\end{equation}
The coefficient $A{(f,\pi)}$ is rational, and for simple $P$ is
essentially a rank $l$ Dedekind sum.  Our next result shows that if
$P$ is nonsingular, then these individual terms transform under the
Hecke operators exactly as the coefficients of $E (P)$ do:

\begin{theorem}\label{thm:main2}
Let $P\in \PPP (L)$ be a nonsingular lattice polytope.  For any
superlattice $M\supset L$, let $f_{M}$ be the face $f$, thought of as
a face of $P_{M}$.  Then for each degree $l$ term $A{(f,\pi
)}\partial_{f}^{\pi }\in \Td_{l} (\Sigma,\dee{h})$ in the Brion-Vergne
formula, we have
\begin{equation}\label{eq:term.eq}
\sum_{M\in \LLL_{k}} A{(f_{M},\pi )}\partial_{f_{M}}^{\pi }\Vol P_{M}
(h)\bigr|_{h=0} = \nu_{n,k,n-l} (p) A{(f,\pi)}\partial_{f}^{\pi }\Vol
P (h)\bigr|_{h=0}.
\end{equation}
\end{theorem}
Note that the Hecke images $P_{M}$ in \eqref{eq:term.eq} are in
general singular, even if $P$ is nonsingular.  Also, the proof of
Theorem \ref{thm:main2} is independent from that of Theorem
\ref{thm:main1}, and hence provides another proof Theorem
\ref{thm:main1} for nonsingular lattice polytopes.

\subsection{} We comment briefly on the proofs of
Theorems~\ref{thm:main1} and~\ref{thm:main2}.  The proof of Theorem
\ref{thm:main1} is a counting argument.  The new lattice points
appearing in $P$ in the sum \eqref{eq:heckeact} all
lie in the superlattice $p^{-1}L$, and to compute $T (p,k)E (P)$ one
keeps track of which lattice points appear in
a given Hecke image.  This gives an expression for $T (p,k)E (P)$ in
terms of $E (P) (t)$, $E (P) (pt)$, and the cardinalities of some
finite Grassmannians.  An additional argument shows that this
expression implies \eqref{eq:thm}.

The proof of Theorem \ref{thm:main2} is more complicated.  At the
heart of \eqref{eq:term.eq} are certain ``distribution relations'' of
Dedekind sums, essentially coming from a distribution relation
satisfied by the Hurwitz zeta function (\S\ref{ss:defofcircle}).  In
the proof of Theorem \ref{thm:main2}, these relations appear in
identities involving Dedekind sums and the cardinalities of strata in
certain stratifications of finite Grassmannians.

Rather than proving these identities directly, we show that they occur
in the computation of the constant term of $T (p,j)E (P')$ for
lower-dimensional polytopes $P'$ and for $j\leq k$.  Since these
constant terms are always $1$, by appropriately choosing $P'$ we show
that our identities hold.  Then we use induction to complete the
argument.

\subsection{}
Here is a fanciful interpretation of Theorem \ref{thm:main1}.  The
Ehrhart polynomial is clearly invariant under the action of $\GL (L)$,
the linear automorphisms of $V$ preserving $L$.  One can think of $\PPP (L)$ as being
like the upper halfplane $\fH$, and the equivalence class of $P\in
\PPP (L)$ as being a point on the modular curve $\SL_{2} (\ZZ)
\backslash \fH$. Then the $l$th coefficient $c_{l}$, thought of as a
function $\GL (L)\backslash \PPP (L) \rightarrow \QQ$, plays the role
of a weight $l$ modular form, and Theorem \ref{thm:main1} says that
$c_{l}$ is a ``weight $l$ Hecke eigenform of level $1$.''
Furthermore, the simple description of its Hecke eigenvalues indicates
that $c_{l}$ should be thought of as being like an Eisenstein series.
Indeed, the analogy between coefficients of $E$ and modular forms was
our original motivation to consider this problem.  These
reflections lead to natural questions unanswered in this paper: 
\begin{itemize}
\item What is the dimension of the space of eigenforms?  Is it
finite-dimensional?
\item What are the analogues of level $N$ modular forms?
\item Are there analogues of modular forms over number fields of
higher degree, e.g. Hilbert modular forms? 
\end{itemize}

\subsection{}
The paper is organized as follows.  Section \ref{s:hecke} recalls
background about lattice polytopes and their normal fans, and
discusses the connection between Hecke operators and finite
Grassmannians.  Section \ref{s:thm1} gives the proof of Theorem
\ref{thm:main1}.  Section \ref{todd.operator} discusses the
computation of the Ehrhart polynomial using the Todd operator, and
Section \ref{s:thm2pf} gives the proof of Theorem \ref{thm:main2}.
Section \ref{s:app} discusses explicit examples of Theorem
\ref{thm:main2} for three-dimensional polytopes, and relates the
identities occurring in the proof of Theorem \ref{thm:main2} to
Dedekind sums and the Hurwitz zeta function.  Finally, Section
\ref{s:avg} addresses the problem of computing the average Ehrhart
polynomial as one varies over a family of superlattices.

\subsection{Acknowledgments} Initial discussions about this project
took place at the Banff International Research Station (BIRS), at the
2003 program \emph{The many aspects of Mahler's measure}.  It is a
pleasure to thank BIRS for its hospitality.  We also thank Michel
Brion for the proof of Lemma \ref{lem:brion.lemma}, and Noam Elkies
and Tom Braden for helpful comments.


\section{Hecke operators and finite Grassmannians}\label{s:hecke}

\subsection{}
Let $P$ be a simple
lattice polytope in the vector space $V$ with vertices in the lattice
$L$.  For convenience we fix a nondegenerative bilinear form $\langle
\phantom{a}, \phantom{b}\rangle$ and use it to identify $V$ with its
dual.  We also assume that $L$ is self-dual with respect to this form.

Let $\FFF$ be the set of faces of $P$, and for any $l$ let $\FFF (l)$
be the subset of faces of dimension $l$.  Let $F\in \FFF (n-1)$ be a
facet of $P$.  Then $F$ is the intersection of $P$ with an affine
hyperplane 
\[
H_{F} = \{x\mid \langle x,u_{F}\rangle +\lambda_{F} =0\},
\]
where the normal vector $u_{F}$ is taken to be a primitive vector
in $L$, and points into the interior of $P$.

\subsection{}\label{ss:nor.cones}
Let $f\in \FFF (n-l)$ be a face of codimension $l$, and let $H_{f}$ be
the affine subspace spanned by $f$.  Since $P$ is simple, there are
exactly $l$ hyperplanes in $\{H_{F}\mid F\in \FFF (n-1) \}$ whose
intersection is $H_{f}$.  Let $\sigma_{f}\subset V$ be the convex cone
generated by the corresponding normal vectors $\{u_{F} \}$.  The cone
$\sigma_{f}$ is called the \emph{normal cone} to $f$.

The set $\Sigma$ of all normal cones $\{\sigma_{f}\mid f\in \FFF
\}$ forms an acute rational polyhedral fan in $V$.  This means the following:
\begin{enumerate}
\item Each $\sigma \in \Sigma$ contains no nontrivial linear subspace.
\item If $\sigma'$ is a face of $\sigma \in \Sigma$, then
$\sigma '\in \Sigma$.  
\item If $\sigma$, $\sigma '\in \Sigma$, then $\sigma \cap \sigma '$
is a face of each.
\item Given $\sigma \in \Sigma$, there exists a finite set $S\subset
L$ such that any point in $\sigma$ can be written as $\sum
\rho_{s}s$, where $s\in S$ and $\rho_{s}\geq 0$.  
\end{enumerate}
Moreover, $P$ simple implies $\Sigma$ is simplicial, which means that
in (4) we can take $\# S = \Dim \sigma$ for all $\sigma$.  The fan
$\Sigma$ is called the \emph{normal fan} to $P$.
 
\subsection{} Let $\rho \in \Sigma$ be a $1$-dimensional cone.  Then
$\rho$ contains a unique normal vector $u_{F}$, which we call
the \emph{spanning point} of $\rho$.  For any cone $\sigma $, we
denote by $\sigma (1)$ the set of spanning points of all
$1$-dimensional faces of $\sigma$, and write 
\[
\Sigma (1) =
\bigcup_{\sigma \in \Sigma} \sigma (1).
\]
There is bijection between
$\Sigma (1)$ and $\FFF (n-1)$.

For any rational cone $\sigma$, let $U (\sigma)$ be the sublattice of
$L$ generated by the spanning points of $\sigma$.  Put $L
(\sigma) = L\cap (U (\sigma)\otimes \QQ )$, and let $\Ind \sigma =
[L(\sigma )\colon U (\sigma)]$.  If $\Ind \sigma = 1$, then
$\sigma$ is called \emph{unimodular}.  Then $P$ is 
nonsingular if and only if all its normal cones are
unimodular.  

\subsection{} Now we recall some basic facts about Hecke operators for
the linear group $\GL_{n}$.  Let $p$ be a prime, and let $\overline V$
be the finite vector space $\FF_{p}^{n}$.  For any rational subspace
$W\subset V$, let $\overline W$ be the corresponding subspace of
$\overline V$.  Fix a positive integer $k\leq n$, and let $\Gr (k,n)$
be the Grassmannian of $k$-dimensional subspaces of an $n$-dimensional
vector space.

\begin{lemma}\label{lem:cosets}
The set $\LLL_{k}$ of superlattices $p^{-1}L\supsetneq M \supsetneq L$
of coindex $p^{k}$ is in bijection with the set $\TTT $ of upper
triangular matrices of the form
\[
\left(\begin{array}{cccc}
p^{e_{1}}&&a_{ij}\\
&\ddots&&\\
&&p^{e_{n}}&\\
\end{array} \right),
\]
where
\begin{itemize}
\item $e_{i}\in \{0,1 \}$, and exactly $k$ of the $e_{i}$ are equal to
$0$, and
\item $a_{ij}=0$ unless $e_{i}=0$ and $e_{j}=1$, in which case
$a_{ij}$ satisfies $0\leq a_{ij}<p$.
\end{itemize}
Moreover, the map $M \mapsto \ol{M}$ induces a bijection between $\LLL_{k}$
and $\Gr (k,n) (\FF_{p})$.
\end{lemma}

\begin{proof}
It is well known that the set of \emph{sublattices} $L\supsetneq N
\supsetneq pL$ of index $p^{n-k}$ is in bijection with $\TTT$
\cite[Prop. 7.2]{krieg}.  To realize this bijection, we take
$L=\ZZ^{n}$, and then any $N$ is constructed as the sublattice
generated by the rows of some $A\in \TTT$.  The sublattice $N$
determines a subspace $\overline N\subset \overline V$, which is the
subspace generated by the $k$ rows with diagonal entry $1$.  It is
clear that we obtain all $k$-dimensional subspaces of $\overline V$ in
this way, for example by considering the decomposition of $\Gr (k,n)
(\FF_p)$ into Schubert cells \cite[p. 147]{fulton}.  Finally, both
statements of the lemma follow from the isomorphism $p^{-1}L/L \simeq
L/pL$ given by scaling by $p$, and from the fact that a sublattice has
coindex $p^{k}$ if and only if it has index $p^{n-k}$.
\end{proof}

\subsection{} Let $f\in \FFF$ be a face of $P$, and let $\sigma_{f}$
be the normal cone to $f$.  Let $V_{f}\subset V$ be the linear
subspace parallel to $H_{f}$, and let $C_{f}$ be the linear span of
$\sigma_{f}$.  The subspace $C_{f}$ contains the distinguished
$1$-dimensional subspaces $\{C_{\rho}\mid \rho \in \sigma_{f} (1) \}$.

\begin{proposition}\label{prop:A}
Let $M\in \LLL_{k}$, and for any $f\in \FFF$, let $f_{M}$ be the corresponding
face of $P_{M}$.  Then
\begin{enumerate}
\item $\Vol f_{M} = p^{\Dim (\overline M \cap \overline V_{f})} \Vol f$, and 
\item $\Ind \sigma_{f_{M}} = p^{\Dim (\overline M \cap \overline C_{f}) - r} \Ind \sigma $,
\end{enumerate}
where 
\[
r = \# \{\overline C_{\rho}\mid \text{$\rho \in \sigma_{f} (1) $ and $\overline
C_{\rho} \subset \overline M$} \}.
\]
\end{proposition}

\begin{proof}
Choose a $\ZZ$-basis $B$ of $L$ such that $B\cap V_{f}$ is a
$\ZZ$-basis for $L\cap V_{f}$.  By Lemma
\ref{lem:cosets}, with respect to $B$ any $M\in \LLL_{k}$ is spanned by the rows of
$p^{-1}A$ for some $A\in \TTT$.  Each row of $A$ with diagonal
entry $1$ contributes a factor of $p$ to $\Vol f_{M}/ \Vol f$, which
proves (1).

For $C_{f}$ we argue similarly.  The only difference is that each row
of $A$ with diagonal entry $1$ contributes a factor of $p$ to $\Ind
\sigma_{f_{M}}/\Ind \sigma_{f}$, unless the diagonal entry is the only
nonzero entry in the row.  This situation corresponds to some subspace
$\overline C_{\rho}$ being contained in $\overline M$, and (2)
follows.
\end{proof}

Proposition \ref{prop:A} allows us to give a geometric interpretation
for the eigenvalue $\nu (p)$.

\begin{proposition}\label{prop:vol.tfm}
Fix nonnegative integers $k,l\leq n$, and let $p$ be a prime. 
Let $\Vol_{l}\colon \PPP (L) \rightarrow \QQ$ be the function 
\[
\Vol_{l}(P) = \sum_{f\in \FFF (l)} \Vol (f),
\]
and define 
\[
T (p,k)\Vol_{l} (P) = \sum_{{M\in \LLL_{k}}} \Vol_{l} (P_{M}).
\]
Then $T (p,k) \Vol_{l} (P) = \nu_{n,k,l} (p) \Vol_{l} (P)$.
\end{proposition}

\begin{proof}
Suppose $f\in \FFF (l)$.
According to Proposition \ref{prop:A}, we have 
\begin{equation}\label{eq:vol.eq}
\sum_{M\in \LLL_{k}} \Vol f_{M} = \sum_{M\in \LLL_{k}} p^{\Dim (\overline M \cap \overline V_{f})} \Vol f.
\end{equation}
The right of \eqref{eq:vol.eq} equals $\nu_{n,k,l} (p)\Vol f$, and the
statement follows immediately.
\end{proof}


\section{Proof of Theorem \ref{thm:main1}}\label{s:thm1}

\subsection{}
Throughout this section we allow $P$ to be
a general lattice polytope.  Let $U\subset \ol{V}$ be a fixed
subspace of dimension $l$ as in \S\ref{ss:nu.def},
and recall 
\[
\nu_{n,k,l}(p) = \sum_{\substack{W\subset \FF_{p}^{n}\\\Dim W=k}}
p^{\Dim W\cap U}. 
\]
Let $G_{k,n}$ be the cardinality of number the finite Grassmannian
$\Gr (k,n) (\FF_p)$.  It is
well known that
\begin{equation}\label{zero}
G_{k,n} = \frac{\pint{n}!}{\pint{k}!\pint{n-k}!},
\end{equation}
where $\pint{n}= (p^{n}-1)/ (p-1)$, and $\pint{n}!=\prod_{i=1}^{n}\pint{i}$.

\begin{lemma}\label{lem:count.arg}
Let $E=E (t)$ be the Ehrhart polynomial of $P$.  Then
\begin{equation}\label{eq:rel1}
T (p,k)E (t) = G_{k-1,n-1}E (pt) + (G_{n,k}-G_{k-1,n-1})E (t).
\end{equation}
In particular,
\begin{equation}\label{eq:rel2}
c_{l} (T (p,k)E)/c_{l} (E ) = G_{k,n}+ (p^{l}-1) G_{k-1,n-1}.
\end{equation}
\end{lemma}

\begin{proof}
We have 
\begin{equation}\label{eq:union}
\bigcup_{M\in \LLL_{k}} M = p^{-1}L,
\end{equation}
and since counting points in $p^{-1}L\cap P$ is done by $E (pt)$,
we must count how often a point $x\in p^{-1}L$ appears in the
union \eqref{eq:union}.  There are two separate cases, namely (i)
$x\in p^{-1}L\smallsetminus L$, and (ii) $x\in L$.  The former
contribute to $E (pt)$, and the latter to $E (t)$.

For (i), note that the point $x$ determines a line $\Lambda_{x}\in
\ol{V}$, and the number of $k$-dimensional subspaces containing
$\Lambda_{x}$ is $G_{k-1,n-1}$.  For (ii), each $x\in L$ will appear in
every Hecke image, which gives $G_{k,n}$ in total.  However, such
points are also counted in the sublattices contributing to (i).  When
these contributions are subtracted, we obtain \eqref{eq:rel1}.  This
proves the first statement.

Finally, \eqref{eq:rel2} follows easily from \eqref{eq:rel1},
since $c_{l} (E (pt)) = p^{l}c_{l}(E (t))$. 
\end{proof}

\begin{lemma}\label{lem:grass.strat}
We have 
\begin{equation}\label{one}
\nu_{n,k,l}(p) = G_{k,n}+ (p^{l}-1) G_{k-1,n-1}.
\end{equation}
Moreover, 
\begin{equation}\label{onea}
\nu_{n,k,l}(p)/\nu_{n,n-k,n-l}(p) = p^{k+l-n}.
\end{equation}
\end{lemma}

\begin{proof}
We treat the case $k\geq l$; the case $k<l$ is similar.

For $j=0,\dotsc ,l$, let $Y_{j}$ be the locally closed subvariety of $\Gr(k,n) (\FF_p)$
defined by
\[
Y_{j} = \{W\mid \Dim W = k, \quad \Dim (W\cap U) =j \},
\]
and let $y_{j}=\#Y_{j}$.  Note that $\sum_{j\geq 0} y_{j} =
G_{k,n}$, and that $\nu_{n,k,l}(p) = \sum_{j\geq 0} y_{j}p^{j}$.
Since $y_{0} = G_{k,n}-\sum_{j\geq 1} y_{j}$, it follows that
\begin{equation}\label{two}
\nu_{n,k,l}(p) = G_{k,n}+ \sum_{j\geq 1} y_{j}(p^{j}-1). 
\end{equation}
We prove the lemma by showing 
\begin{equation}\label{three}
\pint{l}G_{k-1,n-1} = \sum_{j\geq 1} \pint{j}y_{j},
\end{equation}
which is equivalent to \eqref{one} and \eqref{two} taken together.  To
do this, we explicitly describe $Y_{j}$ recursively in terms of
$\{Y_{i}\mid i>j \}$, and show that the right of \eqref{three}
telescopes to the left of \eqref{three}.

Consider first $Y_{l}$.  Any point in $Y_{l}$ is given by choosing a
$k$-dimensional subspace $W$ in $\ol{V}$ containing $U$.  Such
subspaces are in bijection with $(k-l)$-dimensional subspaces of
$\ol{V}/U$, and thus $y_{l} = G_{k-l,n-l}$.

Next, any point in $Y_{l-1}$ is given by choosing an $(l-1)$-dimensional
subspace $S$ of $U$, and then choosing a $k$-dimensional subspace $W$
of $\ol{V}$ with $W\cap U = S$.  The subvariety of those $W$ with
$W\cap U \supset S$ gives $G_{l-1,l} G_{k- (l-1), n- (l-1)}$ points;
this is not $y_{l-1}$ since for each $S$ we have included those $W$
that contain $U$, instead of just meeting $U$ in a subspace of codimension
$1$.  The correct value of $y_{l-1}$ is given by subtracting the
contributions corresponding to points in $Y_{l}$, which gives
\[
y_{l-1}
= G_{l-1,l} (G_{k- (l-1),n- (l-1)} - G_{k-l,n-l}).
\]

For the general $Y_j$ similar considerations apply.  We summarize the
results as follows.  For $j=1,\dotsc ,l$ let $U_{j}\subset
\FF_{p}^{n-j}$ be a fixed subspace of dimension $l-j$, and let $Z_{l}$
be the subvariety of the Grassmanian $\Gr (k-j,n-j) (\FF_p)$ of all
$(k-j)$-dimensional subspaces $W$ such that $W\cap U_{j} = \{0 \}$.
Putting $z_{j} = \#Z_{j}$, we have
\[
z_j =
\begin{cases}G_{k-l,n-l}&j=l,\\
              G_{k-j,n-j}-\sum_{i=1}^{l-j}G_{i,l-j}z_{i+j}&j<l. 
\end{cases}              
\]
Then
\[
y_{j}=G_{j,l}z_{j},\quad \text{$j=1,\dotsc ,l$},
\] 
and in particular
\begin{equation}\label{four}
y_{1}=G_{1,l} (G_{k-1,n-1}-G_{1,l-1}z_{2}-G_{2,l-1}z_{3}-\dotsb
-G_{l-1,l-1}z_{l}).
\end{equation}
Finally, using
\eqref{zero} we see 
\begin{equation}\label{eq:x}
\pint{1}G_{1,l}G_{k-1,n-1}=\pint{l}G_{k-1,n-1},
\end{equation}
and
\begin{equation}\label{eq:y}
\pint{j}G_{j,l} = \pint{1}G_{1,l}G_{j-1,l-1}.
\end{equation}
Using \eqref{eq:x} and \eqref{eq:y} with \eqref{four} shows that the right of
\eqref{three} telescopes to the left of \eqref{three}, which proves \eqref{one}.
A simple computation obtains \eqref{onea} from \eqref{one}, 
and Lemma
\ref{lem:grass.strat} is proved.
\end{proof}

Lemmas \ref{lem:count.arg} and \ref{lem:grass.strat} imply almost all
of Theorem \ref{thm:main1}.  Equations \eqref{eq:rel2} and \eqref{one}
imply \eqref{eq:thm}, and the existence of the polynomial
$\Phi_{n,k,l}$ from \eqref{eq:poly.function} is clear from
\eqref{zero} and \eqref{one}.  The only remaining statement is the
positivity of the coefficients of $\Phi_{n,k,l}$.  To see this, fix a
complete flag in $\ol{V}$
\[
\{0 \}=U_{0}\subsetneq U_{1} \subsetneq \dotsb \subsetneq U_{n} = \ol{V}, 
\]
where $\Dim U_{j} = j$.  We define a polynomial $\widehat{\Phi}_{n,k} \in
\ZZ [x_{0},\dotsc ,x_{n}]$ by 
\begin{equation}\label{eq:phihat}
\widehat{\Phi}_{n,k} = \sum_{\substack{W\subset \ol{V}\\\Dim W=k}}
\prod x_{j}^{\Dim W\cap U_{j}}.
\end{equation}
Clearly $\Phi_{n,k,l} (p)$ is obtained from $\widehat{\Phi}_{n,k}$ by
the substitutions $x_{l}=p$ and $x_{j}=1$ if $j\not =l$.  We claim
$\widehat{\Phi}_{n,k}$ is a polynomial with positive coefficients.
Indeed, the distinct monomials $x^{\alpha} := \prod
x_{j}^{\alpha_{j}}$ in \eqref{eq:phihat} correspond to the different
possibilities of intersections of $W$ with the fixed flag, which
correspond to the decomposition of $\Gr (n,k) (\FF_{p})$ into Schubert
cells $S_{\alpha }$ \cite{fulton}.  Thus we can rewrite \eqref{eq:phihat} as
\[
\widehat{\Phi}_{n,k} = \sum_{\alpha } \#S_{\alpha} (\FF_p) x^{\alpha }.
\]
But each Schubert cell is isomorphic to an affine space, and hence the
coefficients $\#S_{\alpha} (\FF_p)$ are pure $p$-powers.  This
completes the proof of Theorem \ref{thm:main1}.

\begin{remark}\label{rem:tom}
We have the following additional geometric interpretation of the
eigenvalue $\nu_{n,k,l} (p)$.  Let $T$ be the total space of the rank
$n$ trivial bundle over $G (k,n) (\FF_{p})$, and let $T_{l}\subset T$
be the subbundle corresponding to a fixed $l$-dimensional subspace.
Let $B$ be the total space of the tautological bundle over $G (k,n)
(\FF_{p})$, i.e. for any $x\in G (k,n) (\FF_{p})$ the fiber $B_{x}$
over $x$ is the $k$-dimensional subspace corresponding to $x$.  Then
\[
v_{n,k,l} (p) = \# (B \cap T_{l}).
\]
\end{remark}


\section{The Todd operator}\label{todd.operator}

\subsection{}\label{ss:fourone} In this section we describe the Todd operator $\Td
(\Sigma ,\dee{h})$ and how it can be used to compute the Ehrhart
polynomial of a simple lattice polytope $P$.  We closely follow
\cite{bv1}.

Recall that $\FFF$ is the set of faces of $P$, and that each facet
$F\in \FFF (n-1)$ determines an affine hyperplane
\[
H_{F} = \{x\mid \langle x,u_{F}\rangle+\lambda_{F} =0
\}, 
\]
where the normal vector $u_{F}\in L$ is a primitive vector pointing
into the interior of $P$.

Let $h= (h_{F})_{F\in \FFF (n-1)}$ be a real multivariable indexed by
the facets of $P$, and let $P (h)$ be the convex region determined by
the inequalities
\begin{equation}\label{eq:ineq}
\{\langle x,u_{F}\rangle+\lambda_{F}+h_{F}\geq 0 \mid F\in \FFF (n-1)\}.
\end{equation}
Note that $P (0) = P$.  Then $P (h)$ is isomorphic to $P$ for small
$h$, and thus for small $h$ one can consider the volume $\Vol P (h)$.
The following examples will play an important role in the proof of
Theorem \ref{thm:main2}.

\begin{example}\label{ex:simplex}
Let $e_{1}, \dotsc ,e_{n}$ be the canonical basis of $\RR^{n}$, and
let $e_{0}= 0$.  Let $P=\Delta_{n}$ be the convex hull of the vectors
$\{e_{0},\dotsc ,e_{n} \}$.  Then $\Delta_{n}$ is the $n$-dimensional
simplex.  Let $h_{i}$ be the parameter attached to the facet obtained
by deleting the vertex $e_{i}$.  It is easy to check that
\[
\Vol \Delta_{n} (h) = \bigl(1+\sum_{i=0}^{n}h_{i}\bigr)^{n}/n!.
\]
\end{example}

\begin{example}\label{ex:prod}
Let $P$ and $P'$ be two lattice polytopes, and let $h$ and $h'$ be 
multivariables indexed by their facets.  Then 
\[
\Vol (P\times P') (h,h') = \Vol P (h) \Vol P' (h').
\]
In particular, for the unit $n$-cube $P=(\Delta_{1})^{n}$ we obtain
\[
\Vol P (h) = \prod_{i=1}^{n} (1+h_{i}+h_{i}').
\]
\end{example}

\subsection{}
Let $\Sigma$ be the normal fan to $P$.
For any $\sigma\in \Sigma$, define 
\[
Q (\sigma) = \bigl\{\sum_{s\in \sigma (1)} \rho_{s}s \bigm |0\leq
\rho_{s} < 1\bigr\}.
\]
Note that $\Vol Q (\sigma) = \Ind \sigma$, and $Q (\sigma)\cap 
U(\sigma)=\{0 \}$ if and only if $\sigma$ is unimodular.  Put
\[
\Gamma_{\Sigma } = \bigcup_{f\in \FFF} Q (\sigma_{f})\cap L .
\]
We have $\Gamma_{\Sigma}=\{0 \}$ if and only if $P$ is nonsingular.

\subsection{}
For each $F\in \FFF (n-1)$, let $\xi_F\colon V \rightarrow \RR$ be
the unique piecewise-linear continuous function defined by
\begin{itemize}
\item $\xi_{F} (s) = 1$ if $s\in \Sigma (1)$ is the spanning point
corresponding to $F$, 
\item $\xi_{F} (s') =0$ for all other $s'\in \Sigma (1)$, and 
\item $\xi_{F}$ is linear on all the cones of $\Sigma $.
\end{itemize}
Put $a_F (x) = \exp (2\pi i \xi_{F} (x))$ for all $x\in
V $. 

Suppose $g\in \Gamma_{\Sigma}\cap \sigma$.  Then the pair
$(g,\sigma)$ determines a tuple of roots of unity as follows.  If
$s_{1},\dotsc ,s_{l}$ are the spanning points of $\sigma$, and
$F_{1},\dotsc ,F_{l}$ are the corresponding facets, then we can attach
to $(g,\sigma)$ the tuple $(a_{1} (g),\dotsc ,a_{l} (g))$, where we
have written $a_{i}$ for $a_{F_{i}}$.

\subsection{}\label{ss:berndef}
Let $a$ be a complex number and $x$ a real variable.  We define $\Td
(a,\dee{x})$ to be the differential operator given formally by the
power series
\[
\frac{\dee{x}}{1-a\exp (-\dee{x})} = \sum_{k=0}^{\infty} c (a,k)
\left(\frac{\partial }{\partial x}\right)^{k}.
\]
Note that $c (1,k)=B_{k}/k!$, where $B_{k}$ is the $k$th Bernoulli
number.\footnote{With our conventions the Bernoulli numbers are
$B_{1}=1/2$, $B_{2}=1/6$, $B_{4}=-1/30$, \dots , and $B_{2k-1}=0$ for
$k > 1$.  Note that for many authors $B_{1}=-1/2$,
cf. \S\ref{ss:defofcircle}.}  If $a\not =1$, then $c (a,k)$ is a
rational function in $a$ of degree $-1$ closely related to the $k$th
\emph{circle function of Euler} (\S\ref{ss:defofcircle}).  Table
\ref{tab:cir} gives some examples of the $c (a,k)$.

\begin{table}[htb]
\begin{center}
\begin{tabular}{c||l}
$k$&$c (a,k)$\\
\hline\hline
$1$&${-{1}/{(a-1)}}$\\
$2$&${{-a}/{(a^{2} - {2}a + 1)}}$\\
$3$&${-\fractwo{a^{2}+a}{2a^{3} - {6}a^{2} + {6}a - 2}}$\\
$4$&${-\fractwo{a^{3} +{4}a^{2} +a}{6a^{4} - {24}a^{3} + {36}a^{2} - {24}a + 6}}$\\
\end{tabular}
\end{center}
\bigskip
\caption{\label{tab:cir}The coefficients $c (a,k)$.}
\end{table}

\subsection{} \label{ss:fourseven} Now let $h$ be a multivariable with components $h_{F}$
indexed by the facets of $P$.  Let $g\in \Gamma_{\Sigma}$, and define
\[
\Td (g,\dee{h}) = \prod_{F\in \FFF (n-1)} \Td (a_{F} (g),\dee{h_{F}})
\]
and
\begin{equation}\label{eq:todd.op}
\Td (\Sigma ,\dee{h}) = \sum_{g\in \Gamma_{\Sigma}} \Td (g,\dee{h}).
\end{equation}
We have the following theorem, proved by Khovanskii-Pukhlikov if $P$
is nonsingular, and by Brion-Vergne for general simple lattice
polytopes.

\begin{theorem}\label{thm:computingE}
\cite{kp,bv1}
Suppose $P$ is a simple lattice polytope.  
Then the coefficients of the Ehrhart polynomial $E_{P} (t) = \sum_{i=0}^{n}
c_{i}t^{i}$ are given by 
\[
c_{n-l} = \Td_{l} (\Sigma ,\dee{h}) \Vol P (h)\bigr|_{h=0},
\]
where $\Td_{l} (\Sigma ,\dee{h})$ is the degree $l$ part of
$\Td(\Sigma ,\dee{h})$. 
\end{theorem}

For the connection between coefficients of the Todd operator and
higher-dimensional Dedekind sums, we refer to \cite[\S 9]{barv.pommersheim}.


\section{Proof of Theorem \ref{thm:main2}}\label{s:thm2pf}

\subsection{}
We recall some notation from \S\ref{ss:oneseven}.
Let $f\in \FFF$ be a face of codimension $\leq l$, and let $\pi = (\pi
(F))_{F\supset f}$ be an ordered partition of $l$ indexed by the facets
containing $f$.  We 
expand \eqref{eq:todd.op} as a sum over pairs
\[
\Td_{l} (\Sigma,\dee{h}) = \sum_{(f,\pi)} A{(f,\pi)} \partial_{f}^{\pi},
\]
where 
\[
\partial_{f}^{\pi} = \prod_{F\supset f} (\dee{h_{F}})^{\pi (F)}
\]
and
\begin{equation}\label{eq:defofA}
A (f,\pi) = \sum_{g\in \Gamma \cap \sigma_{f}} \prod_{F\supset f} c (a_{F}
(g),\pi (F)).
\end{equation}
Note that if $\sigma_{f}$ is unimodular, then
\begin{equation}\label{eq:nonsing.coeff}
A (f,\pi) = \prod_{F\supset f} \frac{B_{\pi (F)}}{\pi (F)!}.
\end{equation}
 
\subsection{} Now fix a total ordering on (unordered) partitions of
$l$ by using the \emph{lexicographic} order.  In other words, let $\pi
= \{\pi_{1},\dotsc ,\pi_{j} \}$ and $\pi ' = \{\pi '_{1},\dotsc ,\pi
'_{k} \}$ be two partitions of $l$ with parts arranged in
\emph{nonincreasing} order.  Then we have $\pi < \pi '$ if and only if
there exists an index $m$ with $\pi_{i}=\pi_{i}'$ for $i<m$ and
$\pi_{i} < \pi '_{i}$ for $i\geq m$.  For example, if $l=6$, then in
increasing order (and in obvious notation) the partitions are
\[
1^{6}, \,\,21^{4}, \,\,2^{2}1^{2}, \,\,31^{3}, \,\,321, \,\,3^{2},
\,\,41^{2}, \,\,42, \,\,51, \,\,6.
\]

\subsection{} We say the pair $(f,\pi)$ is \emph{squarefree} if $\pi
(F)= 1$ for all $F\supset f$, and we write $\pi = \one$.  We begin
with two lemmas.  Lemma \ref{lem:bv.lem} gives a geometric
interpretation of the squarefree terms, and Lemma
\ref{lem:brion.lemma} allows us to compute nonsquarefree terms using
squarefree terms.  

\begin{lemma}\label{lem:bv.lem}
Let $P$ be simple.  For any face $f\in \FFF$, we have 
\[
\partial_{f}^{\one} \Vol P (h)\bigr|_{h=0} = \frac{\Vol f}{\Ind \sigma_{f}}.
\]
In particular, if $P$ is nonsingular and $f$ has codimension $l$, then 
\[ 
A (f,\one ) \partial_{f}^{\one} \Vol P (h)\bigr|_{h=0} = \frac{\Vol f}{2^{l}}.
\]
\end{lemma}

\begin{proof}
The first statement is Lemma 4.7 in \cite{bv1}.  The second statment
follows from \eqref{eq:nonsing.coeff} since the Bernoulli number $B_{1}$ is $1/2$, and $\Ind
\sigma_{f} =1$ if $P$ is nonsingular.
\end{proof}

\subsection{} The following result is well known to experts, and is
stated (for nonsingular $P$) in \cite[Theorem, p. 795]{kp}.  For the
convenience of the reader we present a proof for $P$ simple.  For
unexplained concepts from toric geometry, we refer to
\cite{fulton.toric}.  What we will need from Lemma
\ref{lem:brion.lemma} is \eqref{eq:fund.rel}.

\begin{lemma}\label{lem:brion.lemma}\cite{b.lett}
Let $X$ be the projective toric variety associated to the simple
lattice polytope $P$.  Then the rational Chow ring $H^{*}(X, \QQ )$ is
isomorphic to the quotient of
\[
\QQ
\left[\dee{h_{F}}\mid F\in \FFF (n-1)\right]
\]
by the ideal $I$ of differential
operators that annihilate the function $\Vol P (h)$.
\end{lemma}

\begin{proof}
The rational Chow ring $H^*(X,\QQ)$ has generators the classes
of the divisors $[D_{F}]$, $F\in \FFF (n-1)$, and the following relations:

\begin{itemize}
\item square-free monomial relations $\prod_{F\in I} [D_F] = 0$
unless the facets in $I$ intersect
transversally along a face of $P$, and 
\item linear relations $\sum_F \langle w, u_F\rangle [D_F] = 0$,
where $w\in L$.
\end{itemize}

But the analogous relations hold for $\QQ [\dee{h_{F}}\mid F\in \FFF
(n-1)]$ applied to $\Vol P (h)$; for example, the linear relations
express invariance of volume under translation. Thus, we obtain a
surjective homomorphism of graded rings
\[
H^*(X,\QQ)   \longrightarrow    \QQ
\left[\dee{h_{F}}\mid F\in \FFF (n-1)\right]/I, \quad [D_F] \longmapsto \dee{h_{F}},
\]
where the $\dee{h_{F}}$ have degree 2. To show its injectivity, it is
enough (by Poincare duality) to show that all intersection
numbers of the form $[D_{F_{1}}]\dotsb [D_{F_{n}}]$ can be read off the
images of the $D_{F}$. But this follows from the formula
\[
\bigl(\sum_F (\lambda_F + h_F) [D_F]\bigr)^n = \Vol P(h),
\]
where the $\lambda_{F}$ come from the inequalities \eqref{eq:ineq} determining $P(h)$.
Indeed, since the $h_F$ are independent variables, any monomial
of degree $n$ in the $[D_F]$ can be expressed in terms of partial
derivatives of $\Vol P (h)$.
\end{proof}

\subsection{}
Let $w\in L$.  Then by Lemma \ref{lem:brion.lemma} the differential operator
\begin{equation}\label{eq:fund.rel}
\sum_{F\in \FFF (n-1)} \langle w, u_{F} \rangle \dee{h_{F}}
\end{equation}
annihilates $\Vol P (h)$.  Hence if $\pi > \one $, by 
repeatedly applying \eqref{eq:fund.rel} we
can write
\begin{equation}\label{eq:der}
\partial_{f}^{\pi} \Vol P (h) = \varepsilon_{f} (\pi) \sum_{f'} \Bigl
( \prod_{w\in W (f')} \langle w, u_{w} \rangle \Bigr)
\partial_{f'}^{\one}\Vol P (h),
\end{equation}
where the quantities in \eqref{eq:der} satisfy
the following:
\begin{itemize}
\item The integer $\varepsilon_{f} (\pi) \in \{\pm 1 \}$ depends only
on the pair $(f,\pi)$;
\item The sum ranges over a finite set of codimension $l$ faces $f'$,
each of which is contained in $f$;
\item For each $f'$, the set $W (f')\subset L\otimes \QQ$ satisfies
\begin{itemize}
\item $\langle w, u\rangle = 1$ for some $u\in \sigma_{f} (1)$, 
\item $\langle w, v\rangle = 0$ for all $v\in \sigma_{f} 
(1)\smallsetminus \{u \}$;
\end{itemize}
\item The $\{u_{w} \} \subset \Sigma (1)$ are such that for each $f'$,
we have 
\[
\sigma_{f'} (1) = \sigma_{f} (1) \cup \{u_{w}\}_{w\in W (f')};
\]
\item The sets $W (f')$ are ordered and 
\[
\langle w', u_{w} \rangle = 0\quad \text{for all $w<w'$.}
\]
\end{itemize}
We choose and fix an expression of the form \eqref{eq:der} for each pair $(f,\pi)$.

\subsection{}
We are now ready to prove Theorem \ref{thm:main2}.  Our goal is to
show 
\begin{equation}\label{eq:nsfstatement}
\sum_{M\in \LLL_{k}} A (f_{M},\pi) \partial_{f}^{\pi} \Vol P_{M} (h) =
\nu_{n,k,n-l} (p) A (f,\pi) \Vol P (h).
\end{equation}

Let $\nu (p) = \nu_{n,k,n-l} (p)$.  Applying \eqref{eq:der} in
\eqref{eq:nsfstatement} and using Lemma \ref{lem:bv.lem}, we see that
it suffices to verify
\begin{multline}\label{eq:one}
\sum_{M\in \LLL_{k}} A (f_{M},\pi) \sum_{f'_{M}} \Bigl (
\prod_{w\in W (f'_{M})} \langle w, u_{w} \rangle \Bigr) \frac{\Vol
f'_{M}}{\Ind \sigma_{f'_{M}}} \\
= \nu (p)A (f,\pi)\sum_{f'} \Bigl (
\prod_{w\in W (f')} \langle w, u_{w} \rangle \Bigr) \Vol (f').
\end{multline}

Since the faces $f'$ appearing in \eqref{eq:der} are independent of
the lattice $M$, we can interchange the sum over $\LLL_{k}$ and the
sum over $f'_{M}$, and focus on a single $f'$.  Furthermore, Lemma
\ref{lem:cosets} implies that the sum over $M$ in \eqref{eq:one} is really a sum
over $\Gr (k,n) (\FF_p)$.  We construct a stratification $\{X_{ij} \}$
of $\Gr (k,n)
(\FF_{p})$ by defining 
\begin{equation}\label{eq:defofXij}
X_{ij} = \{W\subset \overline{V} \mid \Dim W = k, \,\,\Dim W\cap
\overline V_{f} = i, \,\,\Dim W\cap \overline C_{f} = j \},
\end{equation}
and the left of \eqref{eq:one} becomes
\[
\sum_{i,j} \sum_{\ol{M}\in X_{ij}} A (f_{M},\pi) \Bigl (
\prod_{w\in W (f'_{M})} \langle w, u_{w} \rangle \Bigr) \frac{\Vol
f'_{M}}{\Ind \sigma_{f'_{M}}}.
\]

Now let $S_{j} \subset \overline C_{f}$ be a fixed subspace of
dimension $j$, and put
\begin{equation}\label{eq:defofmij}
m_{ij} = \# \{\overline M \in X_{ij} \mid \overline M \supset S_{j}\}.
\end{equation}
The number $m_{ij}$ is independent of the choice of $S_{j}$.  If
$\overline M\in X_{ij}$, then
\begin{equation}\label{eq:vol}
\Vol f'_{M} = p^{i}\Vol f',
\end{equation}
and equation \eqref{eq:one} becomes 
\begin{multline}\label{eq:withXij}
\sum_{i,j} p^{i}m_{ij} \sum_{{\substack{S\subset \overline C_{f'}\\
\Dim S = j}}} A (f_{S},\pi) (\Ind \sigma_{f'_{S}})^{-1} \prod_{w\in W
(f'_{S})} \langle w, u_{w}
\rangle \\
=  \nu (p)A (f,\pi)
\prod_{w\in W (f')} \langle w, u_{w} \rangle,
\end{multline}
where we have written
\begin{equation}\label{eq:afs}
A (f_{S},\pi) = \sum_{g\in \Gamma \cap \sigma_{f_{S}}} \prod_{F\supset
f} c (a_{F} (g),\pi (F)).
\end{equation}
Note that it makes sense to replace the subscript $M$ with $S$ in
\eqref{eq:withXij} and \eqref{eq:afs}, since $\Ind \sigma_{f'_{M}}$
(respectively $\Gamma \cap \sigma_{f_{M}}$) depends only on $S =
\overline M \cap \overline C_{f'}$ (resp., $\overline C_{f}$).  The
notation $W (f'_{S})$ also makes sense, because all points in $W
(f'_{M})$ are multiples of points in $W (f')$ (in fact they differ at
most by a factor of $p$), and which multiples we take depend only on
$S$.

To verify \eqref{eq:withXij}, we show that for each $j$ the
identity 
\begin{equation}\label{eq:lowerrank}
\sum_{{\substack{S\subset \overline C_{f'}\\
\Dim S = j}}} A (f_{S},\pi) \prod_{w\in W (f'_{S})} \langle w, u_{w}
\rangle  (\Ind \sigma_{f'_{S}})^{-1} \\
=   G_{j,l} A (f,\pi)\sum_{f'} \Bigl (
\prod_{w\in W (f')} \langle w, u_{w} \rangle \Bigr)
\end{equation}
holds.  This will complete the proof of the theorem, since
\[
\sum_{i,j}p^{i}m_{ij}G_{j,l} = \nu (p).
\]

We verify \eqref{eq:lowerrank} by induction on the partition order;
the main idea is to show that \eqref{eq:lowerrank} appears in the computation of
the constant term of $T (p,j)E (P)$ for some easily understood
polytope $P$.  Since we know how the constant terms transform under
the Hecke operators, our identity is forced to hold.  In particular, let
\[
P = \prod_{F\supset f} \Delta_{\pi (F)},
\]
where in the product the facets $F$ are ordered so that $\pi$ has
nonincreasing parts.  Using Examples \ref{ex:simplex} and
\ref{ex:prod}, we see that the highest order terms
contributing to $E (P)$ and $T
(p,j)E (P)$ are those of type $(f,\pi)$, where $f$ is a vertex.  Now 
assume that all weight $l$ terms of type $(f,\pi ')$, with $\pi '<\pi$
satisfy \eqref{eq:lowerrank}.  Since the constant term of $T (p,j)E
(P)$ equals $G_{j,l}$, and since each vertex of $P$ contributes
equally to the constant term, this implies \eqref{eq:lowerrank}.

Hence to complete the proof, we must check \eqref{eq:lowerrank} in
the case $\pi = \one$.  In this case we don't need to apply
\eqref{eq:der}, since the terms are already squarefree. 
Using \eqref{eq:nonsing.coeff}, the identity to be proved is
\begin{equation}\label{eq:sqrfreeid}
\sum_{\substack{S\subset \overline C_{f}\\
\Dim S = j}} (\Ind \sigma_{f_{S}})^{-1} A (f,\one) = \frac{G_{j,l}}{2^{l}}.
\end{equation}
To prove \eqref{eq:sqrfreeid}, we let $P= (\Delta_{1})^{l}$ and
consider the action of $T (p,j)$ on the constant term of its Ehrhart
polynomial.  By Example \ref{ex:prod}, we have 
\[
\Vol P (h) = \prod_{i=1}^{l} (1+h_{i}+h_{i}').
\]
We see from applying $\Td_{l}$ to $\Vol P (h)$ that only squarefree
terms contribute to the constant term of $E (P)$, and that this
contribution is the same for all vertices of $P$ (in fact it's
$2^{-l}$).  Moreover, using the matrices given in Lemma
\ref{lem:cosets}, it's easy to see that only squarefree terms
contribute to the contant term of $T (p,j)E (P)$, and that the
contribution for any vertex $f$ is equal to
\begin{equation}\label{eq:eachvert}
\sum_{M\in \LLL_{j}} (\Ind \sigma_{f_{M}})^{-1}A (f_{M},\one ).
\end{equation}
But under $T (p.j)$ the constant term of $E (P)$ is multiplied by
$G_{{j,l}}$, and because the contribution of each vertex is the same,
we have that \eqref{eq:eachvert} equals $G_{j,l}/2^{l}$.  This
completes the proof of \eqref{eq:sqrfreeid}, and the proof of 
Theorem \ref{thm:main2}.

\begin{remark}\label{rem:singularP}
We expect that Theorem \ref{thm:main2} holds if $P$ is replaced by a
general simple lattice polytope, although the argument presented here
doesn't prove this.  In fact, Theorem \ref{thm:main1} suggests that the
analogous result for a general 
lattice polytope should hold, and indeed for the vector partition
functions studied in \cite{bv2}.
\end{remark}

\begin{remark}\label{rem:plugin}
The role of the polytopes $\prod_{F\supset f} \Delta_{\pi (F)}$ in the
proof of Theorem \ref{thm:main2} is very similar to the role of
``basis sequences'' in the theory of characteristic classes and
genera, cf. \cite[p. 79]{hirz}.  This is not a coincidence, since
the machine behind the computation of $c_{l}$ in Theorem
\ref{thm:computingE} is the Hirzebruch-Kawasaki-Riemann-Roch theorem.
\end{remark}


\section{Examples of distribution relations}\label{s:app}

\subsection{}
In this final section, we give examples of the
identites appearing in the proof of Theorem \ref{thm:main2}, and
directly prove them by exhibiting their connection with special values
of the Hurwitz zeta function.  

\subsection{}\label{ss:defofcircle}
Let $u$ be a real number, and let $k$ be a positive integer.  Consider
the special value of the (symmetrized) Hurwitz zeta function 
\[
\zeta (k,u) = \sideset{}{'}\sum_{m\in \ZZ} \frac{1}{(m+u)^{k}}. 
\]
Here the prime next to the summation means to omit the meaningless
term that arises when $u\in \ZZ$.  The series is absolutely convergent
unless $k=1$, in which case we define the value of $\zeta (1,u)$ to
be the limit of the partial sums with $|m|<C$ as $C\rightarrow
\infty$.  Define the \emph{circle functions} $\theta_{k} (u)$ by the
series expansion
\[
\frac{z}{\exp ( z-2\pi i u)-1} = \sum_{k=0}^{\infty} \theta_{k} (u)\frac{z^{k}}{k!}. 
\]
If $u>0$ and $k>1$, then $\theta_{k} (0) = B_{k}$, the $k$th Bernoulli
number as in \S\ref{ss:berndef}.  However note that $c_{1} (0) = -B_{1}$.

By a result of Euler, we have for all $u$
\begin{equation}\label{eq:specval}
\zeta (k,u) = \begin{cases}\displaystyle-\frac{(2\pi i)^{k}}{k!}\theta_{k} (u)&k>1,\\
                   \displaystyle-\frac{(2\pi i)^{k}}{k!}(\theta_{k}(u)+\frac{1}{2})&k=1.
              \end{cases}
\end{equation}

\subsection{}
Now fix a positive integer $n$, and suppose $k>1$.  It is easy to see
that 
\[
\sum_{j=0}^{n-1} \zeta (k,\frac{j}{n}) = n^{k}\zeta (k,0).
\]
Using \eqref{eq:specval}, this becomes 
\begin{equation}\label{eq:thetarel}
\sum_{j=1}^{n-1} \theta_{k} (\frac{j}{n}) = (n^{k}-1)B_{k}.
\end{equation}
Comparing the definition of $c (a,k)$ from \S\ref{ss:berndef} 
yields
\[
c (a,k) = \frac{(-1)^{k}}{k!}\theta_{k} (u),\quad a = \exp (-2\pi iu),
\]
which in \eqref{eq:thetarel} gives
\begin{equation}\label{eq:crel}
\sum_{j=1}^{n-1} c (\omega^{j},k) =\frac{n^{k}-1}{k!}B_{k},\quad k>1.
\end{equation}
Here we have written $\omega = \exp (2\pi i/n)$ and used the fact
that the sum on the left of \eqref{eq:crel} is real.  In fact,
\eqref{eq:crel} remains true if we take $k=1$.

\subsection{}
Let now $P$ be a $3$-dimensional nonsingular lattice polytope;
we investigate the computation of $T (p,1)$ on $c_{1}$.  We focus on
the squarefree case, since no Dedekind sums arise in the nonsquarefree case.

So let $f$ be an edge of $P$.  The key identity \eqref{eq:nsfstatement} becomes
\begin{equation}\label{eq:st2}
\sum_{M\in \LLL_{1}} A (f_{M}, \one)\frac{\Vol f_{M}}{\Ind
\sigma_{f_{M}}} = \frac{p^{2}+2p}{4}\Vol f.
\end{equation}
We break the coefficient $A = A (f_{M},\one)$ into two parts 
\[
A = \nonsing{A}+\sing{A},
\]
where $\nonsing{A}$ corresponds to $g=0$ in \eqref{eq:defofA}, and
$\sing{A}$ corresponds to $g\not =0$.  The latter term appears only if
$\Ind \sigma_{f_{M}} \not = 1$.  Note that $\nonsing{A} = \frac{1}{4}$.

\begin{figure}[htb]
\psfrag{Vf}{$\overline V_{f}$}
\psfrag{Cf}{$\overline C_{f}$}
\psfrag{C1}{$\overline C_{1}$}
\psfrag{C2}{$\overline C_{2}$}
\begin{center}
\includegraphics[scale=0.3]{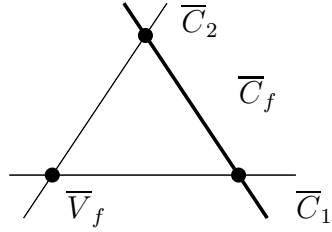}
\end{center}
\caption{\label{fig:ps}Subspaces in $\ol{V}$ for an edge in a $3$-dimensional polytope.}
\end{figure}

To analyze the left of \eqref{eq:st2}, we use Proposition
\ref{prop:A}.  Figure \ref{fig:ps} shows $\overline V$ with the two
subspaces $\overline V_{f}$ and $\overline C_{f}$.  The subspaces
$\overline C_{1}$ and $\overline C_{2}$ are the $1$-dimensional
subspaces corresponding to the two facets containing $f$.  For simplicity, we
draw these subspaces, and the subspaces that follow, by drawing their
images in $\Proj (\overline V) = \Proj^{2} (\FF_p)$.  By abuse of
notation, we denote a subspace of $\overline V$ and the subspace it
induces in $\Proj (\overline V)$ by the same symbol.

\begin{figure}[htb]
\psfrag{1}{$(1,1)$}
\psfrag{2}{$(1,2)$}
\psfrag{3}{$(1,3)$}
\psfrag{4}{$(1,4)$}
\begin{center}
\includegraphics[scale=0.5]{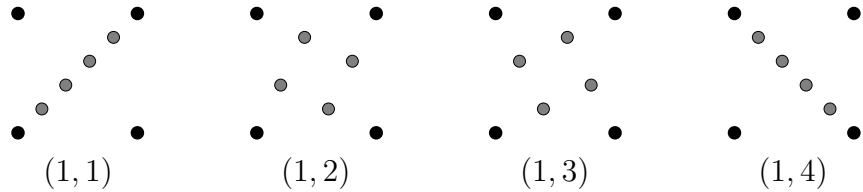}
\end{center}
\caption{\label{fig:t5}Four superlattices giving a nonzero $\sing{A} (\alpha ,\beta)$.}
\end{figure}

Each $M\in \LLL_{1}$ corresponds to a point $\overline M \in \Proj
(\overline V)$.  By Proposition \ref{prop:A}, we have $\Vol f_{M} =
\Vol f$ unless $\ol{M}=\ol{V}_{f}$, in which case $\Vol f_{M} = p \Vol
f$.  Also $\sing{A}=0$ unless $\ol{M}$ meets $\ol{C}_{f}\smallsetminus
\{\ol{C}_{1}\cup \ol{C}_{2} \}$.  Hence there are $p-1$ nonzero
$\sing{A}$, and since $c (a,1)=1/ (1-a)$ each nonzero $\sing{A}$ has
the form 
\[
\sing{A} (\alpha ,\beta) = \sum_{i=1}^{p-1} \frac{1}{(1-\omega^{\alpha
j}) (1-\omega^{\beta j})}, \quad \omega = \exp (2\pi i/p),
\] 
for some nonzero integers $1\leq \alpha ,\beta \leq p-1$.  The value
of $\sing{A} (\alpha ,\beta)$ depends only on the point $[\alpha :
\beta]\in \Proj^{1} (\FF_p)$.  See Figure \ref{fig:t5} for the four
nonzero $\sing{A} (\alpha ,\beta)$ when $p=5$.  The pairs $(\alpha
,\beta)$ are given below each lattice, and the four terms in $\sing{A}
(\alpha ,\beta)$ correspond to the four grey dots.

By \eqref{eq:crel}, the contribution from the singular Hecke images is
\[
\sum_{\substack{[\alpha :\beta]\in \Proj^{1} (\FF_p)\\
[\alpha :\beta] \not = 0, \infty}} \sing{A} (\alpha ,\beta) =
\sum_{i,j=1}^{p-1}\frac{1}{(1-\omega^{i}) (1-\omega^{j})} = \frac{(p-1)^{2}}{4}.
\]
With this in hand it is easy to complete the analysis of
\eqref{eq:st2}.  We break $\Proj (\ol{V})$ into four disjoint subsets 
\[
\Proj (\ol{V}) = S_{1}\cup S_{2}\cup S_{3}\cup S_{4},
\]
where
\begin{itemize}
\item $S_{1} = \ol{V}_{f}$,
\item $S_{2} = \ol{C}_{1}\cup \ol{C}_{2}$,
\item $S_{3} = \ol{C}_{f}\smallsetminus S_{2}$, and 
\item $S_{4} = \Proj (\ol{V}) \smallsetminus \{S_{1} \cup S_{2}\cup S_{3} \}$.
\end{itemize}
The relevant contributions are given in Table \ref{tab:contrib}, and
one easily sees that \eqref{eq:st2} holds.

\subsection{} 
\begin{figure}[htb]
\psfrag{M}{$\ol{M}$}
\psfrag{Vf}{$\overline V_{f}$}
\psfrag{Cf}{$\overline C_{f}$}
\psfrag{C1}{$\overline C_{1}$}
\psfrag{C2}{$\overline C_{2}$}
\begin{center}
\includegraphics[scale=0.3]{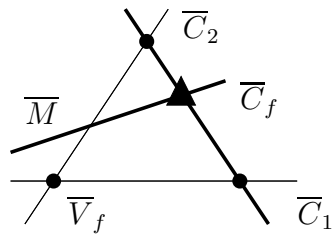}
\end{center}
\caption{\label{fig:lines}Computing $T (p,2)$ on $c_{1}$.}
\end{figure}

The computation of $T (p,2)$ on $c_{1}$ is similar.  The
only difference is that the sum over $M$ corresponds to a sum over
lines in $\Proj (\ol{V})$, and that we obtain a nonzero $\sing{A}$
exactly when a line meets $S_{3}$ in a point.  For example, in Figure
\ref{fig:lines} a nonzero $\sing{A} (\alpha ,\beta )$ arises from the
solid triangle.  Hence each nonzero $\sing{A}(\alpha ,\beta )$ occures
with multiplicity $p$.  Taking this into account, as well as which
lines meet $\ol{V}_{f}$, yields
\[
\sum_{M\in \LLL_{2}} A (f_{M},\one) \frac{\Vol f_{M}}{\Ind
\sigma_{f_{M}}} = \frac{2p^{2}+p}{4}\Vol f.
\]

\begin{table}[htb]
\begin{center}
\begin{tabular}{c||c|c|c|c}
$S_{i}$&$\#S_{i}$&$\Vol f_{M}/\Vol f$&$\Ind
\sigma_{f}/ \Ind \sigma_{f_{M}}$&$\displaystyle\sum_{\ol{M}\in S_{i}} A (f_{M},\one)$\\
\hline\hline
$S_{1}$&$1$&$p$&$1$&$1/4$\\
$S_{2}$&$2$&$1$&$1$&$1/2$\\
$S_{3}$&$p^{2}-1$&$1$&$1$&$(p^{2}-1)/4$\\
$S_{4}$&$p-1$&$1$&$1/p$&$( p^{2}-1 + p-1)/(4p)$
\end{tabular}
\end{center}
\bigskip
\caption{\label{tab:contrib}Summary of $T (p,1)$ on $c_{1}$ for a
$3$-dimensional polytope.}
\end{table}

\section{The regularized Ehrhart polynomial on average}\label{s:avg}

\subsection{}
Let $P$ be a fixed $n$-dimensional lattice polytope respect to the
lattice $L$.  We can define a ``regularized'' version $\tilde{E} (P)$
of $E (P)$ by 
\[
\tilde{E} (P) (t) := E (P) (t) - \Vol (P)t^{n}.
\]

Suppose $\MMM$ is a 
finite set superlattices of $L$ of finite coindex.  We can
define the average regularized Ehrhart polynomial of $P$ with respect to the
family $\MMM$ by 
\[
\tilde{\avg{E}} (P,\MMM) = \frac{1}{\#\MMM}\sum_{M\in \MMM} \tilde{E} (P_{M}). 
\]
Our goal in this section is to show how Theorem \ref{thm:main1} can be
used to derive limiting formulas for
$\tilde{\avg{E}} (P,\MMM)$ as $\MMM$ ranges over families of
superlattices satisfying certain arithmetical conditions.

\subsection{}
As a first example, fix a prime $p$, and suppose $\MMM = \LLL_{1} (p)$ consists of
all superlattices of $L$ of coindex $p$.  Then by definition 
\begin{align*}
\tilde{\avg{E}} (P,\MMM) &= G_{1,n}^{-1}\sum_{l=0}^{n-1} T
(p,1)c_{l}t^{l}\\
&=G_{1,n}^{-1}\sum_{l=0}^{n-1} \nu_{n,1,l} (p) c_{l}t^{l}.
\end{align*}
By Lemma \ref{lem:count.arg}, we have 
\[
\nu_{n,1,l} (p) = G_{1,n}+p^{l}-1 = p^{n-1}+\dotsb +p^{l+1}+2p^{l}+p^{l-1}+\dotsb +p.
\]
This implies the following result:

\begin{proposition}
\[
\lim_{\substack{p\rightarrow \infty\\
p\ \text{\rm prime}}}  \tilde{\avg{E}} (P,\LLL_{1} (p)) =
2c_{n-1}t^{n-1}+c_{n-2}t^{n-2}+\dotsb +c_{1}t + 1.
\]
\end{proposition}

\subsection{}
We can use the relations in the Hecke algebra to derive similar
results for more general sets of superlattices.  Let $T_{p} (n,k)$ be
the operator $T (n,k)$ at the prime $p$, and   
write $T(N) $ for the operator that associates to any lattice $L$ the
set of superlattices of coindex $N$.  Suppose $N$ has prime
factorization $\prod p_{j}^{e_{j}}$. Then, in the algebra $\HHH$
generated by the $T_{p} (n,k)$ as $p$ ranges over all primes $p$, we have 
\cite[Theorem
3.21]{shimura}
\[
T (N) = \prod T (p_{j}^{e_{j}}),
\]
and the operators $T (p^{e})$ satisfy the (formal) identity
\[
\sum_{e=0}^{\infty} T (p^{e})X^{e} = \bigl(\sum_{i=0}^{n} (-1)^{i}p^{i (i-1)/2}T_{p} (n,k)X^{i}\bigr)^{-1}.
\]

As an example of this, suppose $\MMM (p^{2})$ is the set of all
superlattices of $L$ of coindex $p^{2}$.  Note that $\MMM (p^{2}) \not
= \LLL_{2}$, i.e. $T (p^{2})\not = T_{p} (n,2)$.  In fact in $\HHH$ we have the relation 
\[
T (p^{2}) = T_{p} (n,1)^{2}-pT_{p} (n,2).
\]
One can easily show 
\[
\#\MMM (p^{2}) =  G_{1,n}^{2}-pG_{2,n} = G_{2,n+1},
\]
and then from Lemma \ref{lem:count.arg} we find the following:

\begin{proposition}
\[
\lim_{\substack{p\rightarrow \infty\\
p\ \text{\rm prime}}}  \tilde{\avg{E}} (P,\MMM (p^{2})) =
3c_{n-1}t^{n-1}+c_{n-2}t^{n-2}+\dotsb +c_{1}t + 1.
\]
\end{proposition}


\providecommand{\bysame}{\leavevmode\hbox to3em{\hrulefill}\thinspace}
\providecommand{\MR}{\relax\ifhmode\unskip\space\fi MR }
\providecommand{\MRhref}[2]{%
  \href{http://www.ams.org/mathscinet-getitem?mr=#1}{#2}
}
\providecommand{\href}[2]{#2}


\begin{thebibliography}{10}

\bibitem{barv.pommersheim}
A.~Barvinok and J.~E. Pommersheim, \emph{An algorithmic theory of lattice
  points in polyhedra}, New perspectives in algebraic combinatorics (Berkeley,
  CA, 1996--97), Cambridge Univ. Press, Cambridge, 1999, pp.~91--147.

\bibitem{b.lett}
M.~Brion, 2003, personal communication.

\bibitem{bv1}
M.~Brion and M.~Vergne, \emph{Lattice points in simple polytopes}, J. Amer.
  Math. Soc. \textbf{10} (1997), no.~2, 371--392.

\bibitem{bv2}
\bysame, \emph{Residue formulae, vector partition functions and lattice points
  in rational polytopes}, J. Amer. Math. Soc. \textbf{10} (1997), no.~4,
  797--833.

\bibitem{cs}
S.~E. Cappell and J.~L. Shaneson, \emph{Genera of algebraic varieties and
  counting of lattice points}, Bull. Amer. Math. Soc. (N.S.) \textbf{30}
  (1994), no.~1, 62--69.

\bibitem{carlitz}
L.~Carlitz, \emph{A note on generalized {D}edekind sums}, Duke Math. J.
  \textbf{21} (1954), 399--403.

\bibitem{diaz-robins}
R.~Diaz and S.~Robins, \emph{The {E}hrhart polynomial of a lattice polytope},
  Ann. of Math. (2) \textbf{145} (1997), no.~3, 503--518.

\bibitem{ehrhart}
E.~Ehrhart, \emph{Sur un probl\`eme de g\'eom\'etrie diophantienne lin\'eaire.
  {I}. {P}oly\`edres et r\'eseaux}, J. Reine Angew. Math. \textbf{226} (1967),
  1--29.

\bibitem{fulton.toric}
W.~Fulton, \emph{Introduction to toric varieties}, Princeton University Press,
  Princeton, NJ, 1993.

\bibitem{fulton}
\bysame, \emph{Young tableaux}, London Mathematical Society Student Texts,
  vol.~35, Cambridge University Press, Cambridge, 1997, With applications to
  representation theory and geometry.

\bibitem{guillemin}
V.~Guillemin, \emph{Riemann-{R}och for toric orbifolds}, J. Differential Geom.
  \textbf{45} (1997), no.~1, 53--73.

\bibitem{hirz}
F.~Hirzebruch, \emph{Topological methods in algebraic geometry}, Grundlehren
  der mathematischen {W}issenschaften, no. 131, Springer-Verlag, 1978.

\bibitem{kk}
J.-M. Kantor and A.~Khovanskii, \emph{Une application du th\'eor\`eme de
  {R}iemann-{R}och combinatoire au polyn\^ome d'{E}hrhart des polytopes entiers
  de {${\bf R}\sp d$}}, C. R. Acad. Sci. Paris S\'er. I Math. \textbf{317}
  (1993), no.~5, 501--507.

\bibitem{krieg}
A.~Krieg, \emph{Hecke algebras}, Mem. Amer. Math. Soc. \textbf{87} (1990),
  no.~435, x+158.

\bibitem{jppp}
J.~E. Pommersheim, \emph{Toric varieties, lattice points and {D}edekind sums},
  Math. Ann. \textbf{295} (1993), no.~1, 1--24.

\bibitem{kp}
A.~V. Pukhlikov and A.~G. Khovanski\u{\i}, \emph{The {R}iemann-{R}och theorem
  for integrals and sums of quasipolynomials on virtual polytopes}, Algebra i
  Analiz \textbf{4} (1992), no.~4, 188--216.

\bibitem{shimura}
G.~Shimura, \emph{Introduction to the arithmetic theory of automorphic
  functions}, Publications of the Mathematical Society of Japan, vol.~11,
  Princeton University Press, Princeton, NJ, 1994, Reprint of the 1971
  original, Kano Memorial Lectures, 1.

\bibitem{zag-dede}
D.~Zagier, \emph{Higher dimensional {D}edekind sums}, Math. Ann. \textbf{202}
  (1973), 149--172.

\end{thebibliography}
\end{document}